\documentclass[12pt,reqno]{amsart}

\usepackage{amsmath}
\usepackage{amsfonts}
\usepackage{amscd}
\usepackage{amssymb}
\usepackage{amsthm}
\usepackage{mathdots}
\usepackage{mathtools}
\usepackage{bm}

\usepackage[linktocpage=true]{hyperref}
\usepackage{mathrsfs}

\usepackage{microtype}

\usepackage{tikz,tikz-cd, color}

\usepackage{wasysym, stackengine, makebox, graphicx}
\newcommand\isom{\mathrel{\stackon[-0.1ex]{\makebox*{\scalebox{1.08}{\AC}}{=\hfill\llap{=}}}{{\AC}}}}

\usepackage{adjustbox}
\usetikzlibrary{matrix}
\usetikzlibrary{decorations.pathmorphing}

\tikzset{
  symbol/.style={
    draw=none,
    every to/.append style={
      edge node={node [sloped, allow upside down, auto=false]{$#1$}}}
  }
}

\usepackage{stmaryrd}
\usepackage{enumerate}
\usepackage{xcolor}

\usepackage{booktabs}
\usepackage{float}

\mathchardef\mhyphen="2D 

\newtheorem{thm}{Theorem}
\newtheorem*{thm*}{Theorem}

\newtheorem{prop}[thm]{Proposition}
\newtheorem{cor}[thm]{Corollary}
\newtheorem{thm&defn}[thm]{Theorem \& Definition}

\newtheorem{qst}[thm]{Question}

\newtheorem{theorem}[thm]{Theorem}
\newtheorem{lemma}[thm]{Lemma}
\newtheorem{proposition}[thm]{Proposition}

\newtheorem{conjecture}[thm]{Conjecture}

\theoremstyle{definition}
\newtheorem{defn}[thm]{Definition}
\newtheorem{definition}[thm]{Definition}
\newtheorem*{notation*}{Notation}
\newtheorem{example}[thm]{Example}
\newtheorem{remark}[thm]{Remark}

\theoremstyle{remark}

\newtheorem*{rem*}{Remark}

\numberwithin{equation}{section}
\numberwithin{thm}{section}

\usepackage{fullpage}

\newcommand{\bZ}{\mathbb Z}
\newcommand{\bN}{\mathbb N}
\newcommand{\bQ}{\mathbb Q}
\newcommand{\bR}{\mathbb R}
\newcommand{\bC}{\mathbb C}
\newcommand{\bP}{\mathbb{P}}
\newcommand{\bL}{\mathbb{L}}
\newcommand{\bF}{\mathbb{F}}

\newcommand{\cB}{\mathcal{B}}
\newcommand{\cL}{\mathcal{L}}
\newcommand{\cO}{\mathcal{O}}
\newcommand{\cX}{\mathcal{X}}
\newcommand{\cY}{\mathcal{Y}}

\newcommand{\fX}{\mathfrak{X}}
\newcommand{\fz}{\mathfrak{z}}

\DeclareMathOperator{\rZ}{Z}

\newcommand{\sC}{\mathsf{c}}

\newcommand{\tsK}{\widetilde{\mathsf{K}}}
\newcommand{\sP}{\mathsf{p}}

\newcommand{\oM}{\overline{\mathcal{M}}}
\newcommand{\oX}{\overline{X}}

\newcommand{\halpha}{\widehat{\alpha}}
\newcommand{\hell}{\widehat{\ell}}

\newcommand{\tC}{\widetilde{C}}
\newcommand{\tE}{\widetilde{E}}
\newcommand{\tY}{\widetilde{Y}}

\newcommand{\rel}{\mathrm{rel}}
\newcommand{\ev}{\mathrm{ev}}

\newcommand{\vir}{\mathrm{vir}}
\newcommand{\Bl}{\mathrm{Bl}}
\newcommand{\PT}{\mathrm{PT}}

\newcommand{\GW}{\mathrm{GW}}
\newcommand{\vacuum}{|0\rangle}
\newcommand{\tmid}{\,\middle\vert\,}
\newcommand{\tp}{\tau} 

\DeclareMathOperator{\ZPT}{Z_{PT}}

\DeclareMathOperator{\ZGWp}{Z_{GW}^{\prime}}

\DeclareMathOperator{\PTop}{PT}
\DeclareMathOperator{\PTexc}{PT_{exc}}
\DeclareMathOperator{\GWop}{GW^{\prime}}
\DeclareMathOperator{\GWpexc}{GW_{exc}^{\prime}}

\DeclareMathOperator{\Proj}{Proj}
\DeclareMathOperator{\Sym}{Sym}
\DeclareMathOperator{\Aut}{Aut}
\DeclareMathOperator{\NE}{NE}
\DeclareMathOperator{\ch}{ch}
\DeclareMathOperator{\vdim}{vdim}
\DeclareMathOperator{\Img}{Im}

\DeclareMathOperator{\Hilb}{Hilb}

\newcommand{\ignore}[1]{}

\usepackage{accents}
\newcommand{\vect}[1]{\accentset{\rightharpoonup}{#1}}

\title{Gromov--Witten/Pandharipande--Thomas correspondence via conifold transitions}

\author{Yinbang Lin}
\address[]{School of Mathematical Sciences, Key Laboratory of Intelligent Computing and Applications (Ministry of Education), Tongji University, Shanghai 200092, China}
\email{yinbang\textunderscore lin@tongji.edu.cn}

\author{Sz-Sheng Wang}
\address[]{Department of Applied Mathematics, National Yang Ming Chiao Tung University, Hsinchu 30010, Taiwan}
\email{sswangtw@math.nctu.edu.tw}


\begin{document}

\begin{abstract} 
Given a projective conifold transition of smooth projective threefolds from $X$ to $Y$, we show that if the Gromov--Witten/Pandharipande--Thomas descendent correspondence holds for the resolution $Y$, then it also holds for the smoothing $X$ with stationary descendent insertions.
As applications, we show the correspondence in new cases, especially for Fano threefolds.
\end{abstract}

\maketitle

\section{Introduction}\label{intro_sec}

Inspired and motivated by string theory, curve counting on Calabi--Yau threefolds has been one of the central topics in algebraic geometry for decades. 
There are different approaches to this problem. While Gromov--Witten (GW) theory uses stable maps from curves, Donaldson--Thomas (DT) \cite{DT98,Tho00} and Pandharipande--Thomas (PT) \cite{PT09} theories use sheaves with or without extra structure.
All three theories are conjectured to be equivalent.
The {\em correspondence}, namely the equivalence of two theories, was first stated in terms of GW and DT theories by Maulik, Nekrasov, Okounkov and Pandharipande \cite{MNOPI, MNOPII}. 
On the sheaf theoretic side, the DT/PT correspondence has been proven by Bridgeland \cite{Bridgeland11} and Toda \cite{Toda10} for Calabi--Yau threefolds. 
Since the approaches are very different in nature, the GW/DT or GW/PT correspondence is more difficult to be studied.
Considering the works by Bridgeland and Toda, we will focus on GW and PT theories and study the descendent correspondence.
For simplicity, we state here the correspondence conjecture
for Calabi--Yau threefolds \cite[Conj.\ 3.3]{PT09}, which does not require descendent insertions.
\begin{conjecture}\label{MNOP}
Let $M$ be a Calabi--Yau threefold. For a nonzero curve class $\beta$ in $M$, we have the correspondence
\[
      \ZGWp (M; u)_\beta =\ZPT (M; q)_\beta
\]
under the variable change $- q = e^{iu}$.
\end{conjecture}
Here, a Calabi--Yau threefold $M$ is a smooth projective threefold with a trivial canonical line bundle and $H^1 (\cO_M) = 0$. 
The expressions in the conjecture are generating functions of GW and PT invariants in curve class $\beta$ respectively. We will review their definitions in \S\ref{sec:curve} and the descendent correspondence conjecture (Conjecture~\ref{GWPTconj}).

The most important progress regarding GW/PT correspondence is due to Pandharipande and Pixton. They have proven the correspondence in the following cases: 
\begin{enumerate}[(i)]
    \item\label{exaPP:1} smooth projective toric threefolds \cite[Thm.\ 7]{PP14}; 
    \item\label{exaPP:2} Fano or Calabi--Yau complete intersections in products of projective spaces with even cohomology insertions \cite[Thm.\ 1]{PP17}.
\end{enumerate}
Oberdieck \cite{oberdieck2021marked} introduced 
marked relative invariants
, which provide new tools to study arbitrary cohomology insertions. In the 
stationary 
case, i.e., all descendent insertions are even classes of positive degree, Oblomkov, Okounkov and Pandharipande \cite{OOP20} propose an explicit formula for the GW/PT descendent correspondence via vertex operators.

The purpose of the current paper is to prove the correspondence under {\em conifold transitions} (Definition~\ref{defn:con-trans}). They are examples of {\em extremal} (or {\em geometric}) {\em transitions}. An extremal transition is a process of a crepant resolution $Y \to \oX$ followed by a complex smoothing $\oX \rightsquigarrow X$.
We will denote this by $X\nearrow Y$.
It is speculated \cite{Reid87} that (simply connected) Calabi--Yau threefolds can be related via extremal transitions, see \cite{Morrison99, Rossi06} for a survey.

The following is our main result. 
For the precise formulation, see Theorem~\ref{thm:main}.

\begin{theorem}[Theorem~\ref{thm:main}]\label{mainthm}
Let $X \nearrow Y$ be a projective conifold transition of smooth projective threefolds. If $Y$ satisfies the descendent $\GW/\PT$ correspondence (Conjecture \ref{GWPTconj}), then so does $X$ for descendent insertions \eqref{eq:ins-total}. 
\end{theorem}
These are essentially stationary descendent insertions restricted from the total space of the degeneration (cf.\ Remark \ref{rmk_inv}).

To prove the theorem, the key is Proposition \ref{prop:key}. We briefly explain its content and the strategy. See \S \ref{sim_deg_subsec} for the meanings of various symbols.
Following the strategy in \cite{HL12,LLW18,LR01,LY05}, we will degenerate $X$ to $\widetilde{Y}\cup (\sqcup_i Q_i)$ and $Y$ to $\widetilde{Y}\cup (\sqcup_i \tE_i)$, and apply the degeneration formulas \cite{Li02,LW15,MNOPII}. It essentially follows from dimension counts that each $Q_i$ has no contribution to the invariants. 
On the other hand, $\tE_i$ is the toric compactification of a local curve, for which the GW/PT correspondence is known. Thus, we can relate the invariants of $X$ to those of $Y$ via $\widetilde{Y}$.


The main application of the theorem 
will be to establish the correspondence in new cases using the known ones \eqref{exaPP:1} and \eqref{exaPP:2}. It also provides a possible way to prove the conjectural $\GW/\PT$ correspondence with stationary descendent insertions for Fano threefolds (Remark \ref{rmk:allFano3}). We show it 
holds for 44 deformation families of Fano threefolds (Corollary \ref{cor:non-toric}) and a few classes of smoothings of double solids (Corollary \ref{cor:dcover}). More precisely, we have


\begin{cor}[Corollaries \ref{cor:non-toric} \& \ref{cor:dcover}]\label{cor:intr}
The GW/PT correspondence (Conjecture~\ref{GWPTconj}) holds with descendent insertions \eqref{eq:ins-total}, if $X$ is one the following threefolds:
\begin{enumerate}[(a)]
    \item Fano threefolds in Theorem \ref{thm:Galkin}, 
    \item smooth double covers in Proposition \ref{dcover}.
\end{enumerate}
\end{cor}


After posting our preprint on arXiv, we noticed the one by Pardon \cite[Thm.\ 1.6]{pardon2023}. Using very different methods, he proved the GW/PT correspondence {\em without descendents} for projective threefolds with nef anti-canonical bundles. The missing of descendents from his theory was also pointed out, see the paragraph before \S 1.3 Acknowledgements in \cite{pardon2023}. So, his statements, wonderful as they are, do not cover our result, Corollary \ref{cor:intr}. We also point out the anti-canonical bundles of $X$ and $Y$ for the conifold transition in Theorem \ref{mainthm} are not necessarily nef.

The paper is organized as follows. In Section \ref{sec:curve}, we review the (relative) GW and PT-invariants and degeneration formulas. Section \ref{sec:proof} is devoted to the proof of Theorem~\ref{mainthm}. 
In Section \ref{sec:app}, we provide applications.

\begin{notation*}
For a smooth variety $V$, we will denote the integral Mori monoid by $\NE (V)$, namely the set of effective curve classes in $H_2 (V, \bZ) / \mathrm{tors}$. 
If $V$ is complete and $\beta \in \NE (V)$, we set
\[\sC_\beta = \sC_\beta^V \coloneqq \int_\beta c_1 (T_V) .\] 
\end{notation*}

\section{Gromov--Witten and Pandharipande--Thomas theories}\label{sec:curve}

We will briefly review the GW and PT-invariants, their correspondence, 
and the degeneration formulas. We refer the reader to \cite{Li01, Li02, MNOPII, PT09, LW15, Pandharipande18} for details. 

Let $M$ be a smooth projective threefold. Fix a curve class $\beta\in \NE(M)$, integers $r\in \mathbb{Z}_{\geqslant 0}$ and $n,g\in \mathbb{Z}$.
\subsection{GW and PT-invariants}
We review descendent GW and PT-invariants of threefolds and the corresponding invariants relative to a divisor.

\subsubsection{Absolute theories}
First, let $\oM' _{g, r} (M, \beta)$ denote the moduli space of $r$-pointed stable maps
\begin{align*}
    C\to M
\end{align*}
with possibly disconnected domain curves $C$ of (arithmetic) genus $g$ and \emph{no} contracted connected components (cf.\ \cite{Kon95, FP97}). The latter condition requires each connected component of $C$ to represent a nonzero class in $\NE (M)$ and hence $\beta = [C] \neq 0$.
The moduli space $\oM' _{g, r} (M, \beta)$ is equipped with a virtual fundamental class \cite{BehFan97,LiTia98} and its virtual dimension is $\sC_\beta + r$.
Consider the first Chern class of cotangent line bundle $\bL_i$ associated to the $i$-th marked point:
$\psi_i = c_1 (\bL_i) \in H^2 (\oM' _{g, r} (M, \beta), \bQ)$, $i = 1, \cdots, r$.
Let 
\[\ev_i \colon \oM' _{g, r} (M, \beta) \to M, \quad \mbox{for } 1\leqslant i\leqslant r, \]
be the evaluation maps. Given 
$\gamma_1, \cdots, \gamma_r \in H^\ast (M, \bQ)$, define the disconnected descendent GW-invariants by
\begin{align*}
    \langle \tau_{k_1} (\gamma_1) \cdots \tau_{k_r} (\gamma_r) \rangle' _{g, \beta} = \int_{[\oM'_{g, r} (M, \beta)]^\vir} \prod_{i = 1}^r 
    \psi_i^{k_i} \cup \ev_i^\ast (\gamma_i).
\end{align*}
Note that $\oM' _{g, r} (M, \beta)$ is empty for $g$ sufficiently negative. 
The associated partition function
\begin{align*} 
    \ZGWp \left(M; u \tmid  \prod_{i = 1}^r \tau_{k_i} (\gamma_i) \right)_\beta = \sum_{g \in \bZ} \left\langle \prod_{i = 1}^r \tau_{k_i} (\gamma_i) \right\rangle' _{g, \beta} u^{2g - 2} \in \bQ (\!(u)\!)
\end{align*}
is a Laurent series.


Next we consider the moduli space of stable pairs. A {\em stable pair} \begin{align*}
    (F, s\colon \cO_M \to F)
\end{align*} on $M$ consists of a pure $1$-dimensional sheaf $F$ on $M$ and a section $s$ with $0$-dimensional cokernel. 
Given $n \in \bZ$ and nonzero $\beta \in \NE (M)$, let $P_n (M, \beta)$ be the moduli space of stable pairs with $\ch_2 (F) = \beta$ and $\chi (F) = n$. The moduli space $P_n (M, \beta)$ is fine and projective, and it admits a virtual fundamental class of virtual dimension $\sC_\beta$ \cite[Thm. 2.14]{PT09}. 
Let 
\begin{align*}
    \bF \to  P_n (M, \beta)\times M
\end{align*}
be the universal sheaf.
Let $\pi_P$ and $\pi_M$ be the projections from $P_n (M, \beta) \times M$ onto the first and second factors respectively. For $k \in \bZ_{\geqslant 0}$, the $k$-th descendent insertion $\tp_k (\gamma)$ of a class $\gamma \in H^{p} (M, \bQ)$ is defined by
\begin{align*}
    \tp_k (\gamma) (\xi) = \pi_{P \ast} (\pi_M^\ast (\gamma) \cdot \ch_{2 + k} (\bF) \cap \pi_P^\ast (\xi)) \in H_{\bullet - 2k + 2 - p} (P_n (M, \beta), \bQ)
\end{align*}
for every $\xi \in H_{\bullet} (P_n (M, \beta), \bQ)$.
We use the same symbol to denote descendent insertions in GW and PT theories whose meaning should be clear from the context. We will soon see there is a close relation between the two types of insertions via the GW/PT correspondence.

Given $k_i\in \mathbb{Z}_{\geqslant 0}$ and $\gamma_i \in H^\ast (M, \bQ)$ for $i=1,\cdots,r$, the corresponding descendent PT-invariant is
\begin{align*}
    \langle \tp_{k_1} (\gamma_1) \cdots \tp_{k_r} (\gamma_r) \rangle_{n, \beta} = \int_{[P_n (M, \beta)]^\vir} \prod_{i = 1}^r \tp_{k_i} (\gamma_i).
\end{align*}
Note that the moduli space $P_n (M, \beta)$ is empty for $n$ sufficiently negative. The associated partition function
\begin{align*} 
    \ZPT \left(M; q \tmid \prod_{i = 1}^r \tp_{k_i} (\gamma_i) \right)_{\beta} = \sum_{n \in \bZ} \left\langle \prod_{i = 1}^r \tp_{k_i} (\gamma_i) \right\rangle_{n, \beta} q^n \in \bQ (\!( q )\!)
\end{align*}
is a Laurent series as well.
The following conjecture of the rationality of partition function was made in \cite[Conj.\ 1]{PT09v}.

\begin{conjecture}
\label{rat_conj}
The partition function  $\ZPT \left(M; q \mid \prod_{i = 1}^r \tp_{k_i} (\gamma_i) \right)_{\beta}$ is the Laurent expansion of a rational function in $q$.   
\end{conjecture}

\begin{remark}
If $M$ is Calabi--Yau, then the DT/PT correspondence was proved by Toda \cite{Toda10} (see also \cite[\S 1.2]{Toda20}) and Bridgeland \cite{Bridgeland11}. Moreover, we have the rationality property: $\ZPT (M; q)_{\beta} \in \bQ (q)$ which is invariant under $q \leftrightarrow q^{- 1}$. Hence, the variable change in Conjecture \ref{MNOP} is well-defined.
\end{remark}

\subsubsection{Relative theories}\label{rel_theories_subsusec}

Let $D$ be a smooth (not necessarily connected) divisor on $M$. Relative GW and PT theories enumerate curves with specified tangency to the divisor $D$. To impose the boundary conditions along $D$, we need the notion of cohomology weighted partitions (see also \cite[\S 3.1]{MNOPII}).
 
\begin{definition}
Assume that $D$ is connected, and let $\cB$ be a basis of $H^\ast (D, \bQ)$. A {\em cohomology weighted partition} $\eta$ with respect to $\cB$ is a set of pairs 
\begin{align*}
    \{(a_1, \delta_1), \cdots, (a_{\ell}, \delta_{\ell})\}, \quad \mbox{where } \delta_j \in \cB \mbox{ and } a_1 \geqslant \cdots \geqslant a_{\ell} \geqslant 1,
\end{align*}
such that $\vect{\eta} \coloneqq (a_j) \in \bN^{\ell}$ is a partition of 
size $|\eta| \coloneqq \sum a_j$ and length $\ell (\eta) \coloneqq \ell$. 

The automorphism group $\Aut (\eta)$ consists of $\sigma \in \mathfrak{S}_{\ell (\eta)}$ such that $\eta^{\sigma} = \eta$.
\end{definition}

Let $D_1,\cdots, D_k$ be the connected components of $D$ and 
\[\eta_i = \{(a_{ij}, \delta_{ij})\}_j,\ \mbox{for }1\leqslant i\leqslant k, \]
a cohomology weighted partition over $D_i$ with respect to a fixed basis $\cB_i$ of $H^\ast (D_i, \bQ)$. Let $\beta \in \NE (M)$ be a nonzero curve class satisfying $\beta \cdot D_i = |\eta_i| \geqslant 0$ for each $1 \leqslant i \leqslant k$.

In relative GW-theory, the numbers $a_{ij}$ record the multiplicities of intersection with the connected divisor $D_i$ while the cohomology classes $\delta_{ij}$ record where the tangency occurs. More precisely, we consider the moduli space introduced by J.\ Li \cite[Def.\ 4.9]{Li01}
\begin{align*}
    \oM'_{g, r} (M / D, \beta, \vect{\eta}_1, \cdots, \vect{\eta}_k)
\end{align*}
which parametrizes $r$-pointed relative stable maps of (arithmetic) genus $g \in \bZ$ and degree $\beta$ with possibly disconnected domain curves and relative multiplicities determined by $\vect{\eta}_1, \cdots, \vect{\eta}_k$. 
When the cohomology weighted partitions are empty, we omit them from the expression.
An element in the moduli space is a map to the stack of expanded relative pairs. 
As usual, a relative stable map has nonzero degrees on every connected component of its domain. It carries a virtual fundamental class of (complex) dimension 
\begin{align}\label{rel_vdGW}
    \sC_\beta^M + (\ell (\eta_1) - |\eta_1|) + \cdots + (\ell (\eta_k) - |\eta_k|) + r,
\end{align}
see for example \cite[p.160]{LLZ07}.

For $1\leqslant i\leqslant k$ and $1 \leqslant j \leqslant \ell (\eta_i)$, the moduli space has a relative evaluation map
\begin{align*}
    \ev_{D_i, j} \colon \oM'_{g, r} (M / D, \beta, \vect{\eta}_1, \cdots, \vect{\eta}_k) \to D_i,
\end{align*}
which sends a relative stable map to the $j$-th intersection point with the divisor $D_i$ (according to the fixed ordering). By abuse of notation, we write 
\begin{align*}
    \ev_{D_i}^\ast (\delta_{\eta_i}) \coloneqq \prod_{j = 1}^{\ell (\eta_i)} \ev_{D_i,j}^{\ast} (\delta_{i j}).
\end{align*}
Given 
\begin{align}\label{eq:des-ins}
    \gamma_1, \cdots, \gamma_r \in H^\ast (M, \bQ)\mbox{ and } k_i\in \bZ_{\geqslant 0}, \mbox{ for } 1\leqslant i\leqslant r,
\end{align}
the relative descendent GW-invariants \cite[p.240]{Li02} are
\begin{align*}
    &\langle \tau_{k_1} (\gamma_1) \cdots \tau_{k_r}  (\gamma_r) \mid  \eta_1, \cdots, \eta_k \rangle' _{g, \beta} \\ 
    ={}& \frac{1}{\prod_{j = 1}^k |\Aut (\eta_j)|} \int_{[\oM'_{g, r} (M / D, \beta,  \vect{\eta}_1, \cdots, \vect{\eta}_k)]^{\vir}} \prod_{i = 1}^r \left(\psi_i^{k_i} \cup  \ev_i^\ast (\gamma_i) \right) \cup \prod_{j = 1}^k \ev_{D_j}^\ast (\delta_{\eta_j}).
\end{align*}
Then the associated partition function
\begin{align}\label{relZGW}
    \ZGWp \left(M / D ; u \tmid \prod_{i = 1}^r \tau_{k_i} (\gamma_i) \tmid \eta_1, \cdots, \eta_k  \right)_{\beta} 
    = \sum_{g \in \bZ} \left\langle \prod_{i = 1}^r \tau_{k_i} (\gamma_i) \tmid \eta_1, \cdots, \eta_k \right\rangle'_{g, \beta} u^{2 g - 2}
\end{align}
is a Laurent series as before.

In relative PT-theory, we consider the moduli space introduced by Li-Wu \cite{LW15} (cf.\ \cite[\S 3.2]{MNOPII})
\begin{align*}
    P_n (M / D, \beta)
\end{align*}
which parametrizes stable pairs $(F,s)$ relative to $D$, such that $\chi (F) = n \in \bZ$ and $\ch_2 (F) = \beta$. It carries a virtual fundamental class of dimension \cite[Lem.\ 2]{MNOPII}
\begin{align}\label{eq:rel-pt-vdim}
    \vdim P_n(M/D,\beta)=\sC_\beta^M = \int_\beta c_1 (T_M).
\end{align}
For each $1 \leqslant i \leqslant k$, we have the intersection map
\begin{align}\label{relPTintmap}
    \epsilon_i \colon P_n (M/D, \beta) \to \Hilb( D_i, |\eta_i| )
\end{align}
to the Hilbert scheme of $|\eta_i| = \beta \cdot D_i$ points of the connected divisor $D_i$. 

We recall the Nakajima basis for the cohomology of Hilbert schemes of points and refer the reader to \cite{Nakajima} for more details. Fix $d \in \bN$ and let $\eta = \{ (a_j, \delta_j)\}_j$ be a cohomology weighted partition of size $d$ with respect to $\cB_i$. Set 
\begin{align*}
    \fz (\eta) = |\Aut (\eta)| \cdot \prod_{j = 1}^{\ell (\eta)} a_j.
\end{align*}
Following the notation in \cite{Nakajima} and \cite[\S 3.2.2]{MNOPII}, we write
\begin{align*}
    C_\eta = \frac{1}{\fz (\eta)} P_{\delta_1} [a_1] \cdots P_{\delta_{\ell (\eta)}} [a_{\ell(\eta)}] \cdot \mathbf{1} \in H^{\ast} (\Hilb(D_i, d ), \bQ).
\end{align*}
Here $\mathbf{1}$ is the vacuum vector $\vacuum=1\in H^0(\Hilb(D_i,0), \bQ)$.
Then $\{C_\eta\}_{|\eta| = d}$ is the Nakajima basis of $H^{\ast} (\Hilb(D_i,d), \bQ)$.

\begin{definition}
Suppose that the cohomology basis $\cB_i$ of $D_i$ is self dual with respect to the Poincar\'e pairing, i.e., for each $j$, $\delta_j^\vee =\delta_l$ for some $l$. The {\em dual partition} $\eta^\vee$ is the cohomology weighted partition $\{(a_j, \delta_j^\vee)\}_j$ (with respect to $\cB_i$). 
\end{definition}

Note that the Nakajima basis is orthogonal with respect to the Poincar\'e pairing,
\begin{align*}
    \int_{\Hilb(D_i,d )} C_\eta \cup C_\nu = \begin{cases}
        \frac{(- 1)^{d - \ell (\eta)}}{\fz (\eta)}, &\text{if }\nu=\eta^\vee \\
        0, &\text{otherwise.}
    \end{cases} 
\end{align*}

Given \eqref{eq:des-ins}, the relative descendent PT-invariants are
\begin{align*}
    \langle \tp_{k_1} (\gamma_1) \cdots \tp_{k_r} (\gamma_r) \mid \eta_1, \cdots, \eta_k \rangle_{n, \beta} = \int_{[P_n (M / D, \beta)]^{\vir}} \left(\prod_{i = 1}^r \tp_{k_i} (\gamma_i) \right) \cup \prod_{j = 1}^k \epsilon_j^\ast (C_{\eta_j}).
\end{align*}
The associated partition function is 
\begin{align}\label{relZPT}
     \ZPT \left(M / D ; q\tmid \prod_{i = 1}^r \tp_{k_i} (\gamma_i) \tmid \eta_1, \cdots, \eta_k \right)_{\beta} 
     = \sum_{n \in \bZ} \left\langle \prod_{i =1}^r \tp_{k_i} (\gamma_i) \tmid \eta_1, \cdots, \eta_k \right\rangle_{n, \beta} q^n.
\end{align}
Again, if the cohomology weighted partitions $\eta_j$'s are empty, we will omit them from the expression.
The following rationality statement here is parallel to the absolute case \cite{MNOPII, PP14, PP17}.

\begin{conjecture}\label{relrat_conj}
Assume $D$ is connected. The descendent partition function 
\begin{align*}
    \ZPT \left(M / D; q\tmid \prod_{i = 1}^r \tp_{k_i} (\gamma_i) \tmid \eta  \right)_{\beta} \in \bQ (\!( q )\!)
\end{align*}
is the Laurent expansion in $q$ of a rational function.
\end{conjecture}

\subsection{GW/PT correspondence}
Descendent GW and PT-invariants are very different in flavor. The key to relate them is the correspondence matrices found by Pandharipande and Pixton \cite{PP14,PP17}. See also \cite[\S 1.4 \& 5.4]{oberdieck2021marked}. The matrices relating GW and DT-invariants were predicted in \cite[Conj.\ 4]{MNOPII}.
\subsubsection{Absolute version}
Let $\halpha = (\halpha_1, \cdots , \halpha_{\hell})$, with $\halpha_1 \geqslant \cdots \geqslant \halpha_{\hell} \geqslant 1$, be a partition of length $\ell (\halpha) = \hell$ and size $|\halpha| = \sum \halpha_j$. 
Let $\iota_{\Delta} \colon \Delta \to M^{\hell}$ be the inclusion of the small diagonal in the product $M^{\hell}$. For $\gamma \in H^\ast (M, \bQ)$, we write 
\begin{align*}
    \gamma \cdot \Delta \coloneqq \iota_{\Delta \ast} (\gamma) \in H^\ast (M^{\hell}, \bQ)
\end{align*}
and 
\begin{align*}
    \ev_{\{1, \cdots, \hell\}} \coloneqq (\ev_1, \cdots, \ev_{\hell} ) \colon \oM' _{g, \hell} (M, \beta) \to M^{\hell}.
\end{align*} 
The descendent insertion $\tau_{[\halpha]} (\gamma)$ denotes
\begin{align}\label{tau_part}
    \tau_{[\halpha]} (\gamma) \coloneqq \psi_1^{\halpha_1 - 1} \cdots \psi_{\hell}^{\halpha_{\hell} - 1} \cdot \ev_{\{1, \cdots, \hell\}}^\ast (\gamma \cdot \Delta).
\end{align}
Alternatively, let $\{\theta_j\}$ be a basis of $H^\ast (M, \bQ)$. By K\"unneth formula, we have
\begin{align*}
    \gamma \cdot \Delta = \sum_{j_1, \cdots, j_{\hell}} c^{\gamma}_{j_1, \cdots, j_{\hell}} \theta_{j_1} \otimes \cdots \otimes \theta_{j_{\hell}},
\end{align*}
and then we have \cite[(3)]{PP17}
\begin{align*}
    \tau_{[\halpha]} (\gamma) = \sum_{j_1, \cdots, j_{\hell}} c^{\gamma}_{j_1, \cdots, j_{\hell}} \tau_{\halpha_1 - 1}(\theta_{j_1}) \cdots \tau_{\halpha_{\hell} - 1}(\theta_{j_{\hell}}).
\end{align*}

\begin{example}
If $\gamma$ is the class $\sP$ of a point, then 
\[
    \tau_{[\halpha]} (\sP) = \tau_{\halpha_1 - 1} (\sP) \cdots \tau_{\halpha_{\hell} - 1} (\sP).
\]
If $\halpha = (\halpha_1)$, then $\tau_{[\halpha]} (\gamma) = \tau_{\halpha_1 - 1} (\gamma)$.
\end{example}

A universal correspondence matrix $\tsK$ between the descendent insertions in GW and PT theories was constructed in \cite[\S 0.5]{PP14}. The elements
\begin{align*}
    \tsK_{\alpha, \halpha} \in \bQ [i, c_1, c_2, c_3] (\!( u)\!)
\end{align*}
of the matrix are indexed by partitions $\alpha$ and $\halpha$ of positive size and depend on $i = \sqrt{- 1}$ and the formal variables $c_j$ and $u$. By convention the variable $c_j$ has degree $j$. 
The elements of $\tsK$ satisfy the following two properties \cite[(59) and Prop. 24]{PP14}:
\begin{enumerate}[(a)]
    \item
    If $|\alpha| < |\halpha|$, then  $\tsK_{\alpha, \halpha} = 0$.

    \item
    The $u$ coefficients of $\tsK_{\alpha, \halpha} \in \bQ [i, c_1, c_2, c_3] (\!( u)\!)$ are homogeneous in the variables $c_i$ of degree
    \[
        |\alpha| + \ell (\alpha) - |\halpha| - \ell (\halpha) - 3 (\ell (\alpha) - 1).
    \]
\end{enumerate}
By specializing the formal variables $c_i$ to $c_i (T_M)$, the elements of $\tsK$ act by cup product on $H^\ast (M, \bQ)$ with $\bQ[i](\!(u)\!)$-coefficients:
\begin{align*}
    \tsK_{\alpha, \halpha} \colon H^\ast (M, \bQ) \to H^\ast (M, \bQ[i](\!(u)\!))
\end{align*}
for each partitions $\alpha$ and $\halpha$ of positive size.

Let $\alpha = (\alpha_1, \cdots, \alpha_{\ell})$ be a partition and $P$ a partition of $\{1, \cdots, \ell\}$. For each $S \in P$, a subset of $\{1, \cdots, \ell\}$, let $\alpha_S$ be the subpartition consisting of the parts $\alpha_j$ for $j \in S$ and 
\[
    \gamma_S = \prod_{j \in S} \gamma_j.
\]

\begin{definition}[\cite{PP14}]\label{inser_def_abs}
For even cohomology classes $\gamma_j \in H^{2 \ast} (M, \bQ)$, let
    \begin{align*}
      \overline{\tau_{\alpha_1 - 1} (\gamma_1) \cdots \tau_{\alpha_{\ell} - 1} (\gamma_{\ell})} = \sum_{\substack{P \text{ set partitions } \\ \text{of } \{1, \cdots, \ell \}}} \prod_{S \in P}  \sum_{0 < |\halpha| \leqslant |\alpha_S|} \tau_{[\halpha]} \left(\tsK_{\alpha_S, \halpha} \cdot \gamma_S \right).
    \end{align*}
\end{definition}

\begin{remark}\label{sign_rmk}
In general, a sign has to be included in Definition \ref{inser_def_abs} when there is odd cohomology, see \cite[p.2758]{PP14}. But we will focus on even insertions.
\end{remark}

\begin{example}[\cite{PP14}]
    Let $\alpha = (\alpha_1, \cdots, \alpha_{\ell})$ be a partition and $\gamma_j \in H^{2 \ast} (M, \bQ)$ even cohomology classes.
\begin{enumerate}[(a)]
    \item
    We can write the descendent correspondence as
    \begin{align*}
        \overline{\tau_{\alpha_1 - 1} (\gamma_1) \cdots \tau_{\alpha_{\ell} - 1} (\gamma_{\ell})} = (i u)^{\ell(\alpha) - |\alpha|} \tau_{\alpha_1 - 1} (\gamma_1) \cdots \tau_{\alpha_{\ell} - 1} (\gamma_{\ell}) + \cdots,
    \end{align*}
    where the dots stand for terms $\tau_{[\halpha]}(\cdots)$ with $|\halpha| < |\alpha|$.
    
    \item
    For the case $\alpha = (1^\ell)$, we have 
    \begin{align*}
        \overline{\tau_{0} (\gamma_1) \cdots \tau_{0} (\gamma_{\ell})} = \tau_{0} (\gamma_1) \cdots \tau_{0} (\gamma_{\ell}).
    \end{align*}
\end{enumerate}
\end{example}

We are now in a position to state the conjectural GW/PT correspondence.

\begin{conjecture}[\cite{PP14}]\label{GWPTconj}
Let $\alpha = (\alpha_1, \cdots, \alpha_{r})$ be a partition.
For (even) classes $\gamma_j \in H^\ast (M, \bQ)$, $1\leqslant j\leqslant r$, we have
\begin{align*}
    &(- q)^{- \sC_\beta / 2} \ZPT \left(M  ; q\tmid \tp_{\alpha_1 - 1} (\gamma_1)  \cdots \tp_{\alpha_{r} - 1} (\gamma_{r}) \right)_\beta \\
     ={}&(- iu)^{\sC_\beta} \ZGWp \left(M; u\tmid \overline{\tau_{\alpha_1 - 1} (\gamma_1) \cdots \tau_{\alpha_{r} - 1} (\gamma_{r})}  \right)_\beta 
\end{align*}
under the variable change $- q = e^{iu}$.
\end{conjecture}

The variable change is well-defined assuming
Conjecture \ref{rat_conj}.

For the toric case, \cite[Thm.\ 7]{PP14} implies the following theorem by taking the non-equivariant limit.

\begin{theorem}[\cite{PP14}]\label{GWPTtoric}
If $M$ is a smooth projective toric threefold, then it satisfies the $\GW/\PT$ correspondence, i.e., Conjecture \ref{GWPTconj}.
\end{theorem}

We next review the correspondence over the local $\mathbb{P}^1$.
Let $N = \cO_{\bP^1} (- 1)^{\oplus 2}$ and $P$ the projective bundle\footnote{We are following the classical tradition, $P(E) = \Proj (\Sym E^\vee)$.} $P(N \oplus \cO_{\bP^1})$ over $\bP^1$. Let $\tC\subseteq P$ be the subcurve given by the inclusion $\cO_{\mathbb{P}^1}\to N\oplus \cO_{\mathbb{P}^1}$ and $E$ the hyperplane at infinity in $P$ given by $N \to N \oplus \cO_{\bP^1}$. By the Euler sequence, we have
\begin{align}\label{loccur_T}
    c_1 (T_{P}) = c_1 (\cO_{P} (3)) = 3 [E],
\end{align}
and hence $\int_{\tC} c_1 (T_P) = 0$. Because $P$ is a smooth projective toric threefold, the following statement is a special case of \cite[Thm.\ 4]{PP14} (cf.\ \cite[\S3.3]{MNOPI}), which will be used in the proof of Theorem \ref{thm:main}.

\begin{theorem}\label{GWPTloc}
For each $d \in \bN$, we have the correspondence
\[
    \ZPT (P; q)_{d [\tC]} = \ZGWp (P; u)_{d [\tC]}
\]
under the variable change $- q = e^{iu}$.
\end{theorem}

\subsubsection{Relative version}
A relative version of the correspondence matrix was introduced in \cite[\S 0.4 \& \S 1.3]{PP17}.
Let $D$ be a smooth divisor of $M$.
For $s \in \bN$, let $(M / D)^s$ be the moduli space of $s$ ordered (possibly coincident) points in $M$ relative to $D$: 
\begin{align*}
    p_1, \cdots, p_s \in M / D,
\end{align*}
cf.\ \cite[Def.\ 2.5]{oberdieck2021marked}.
Note that it is proper and smooth of dimension $s \dim M  = 3 s$. Let
\begin{align*}
    \Delta_\rel \subseteq (M / D)^s
\end{align*}
be the small diagonal where all the points $p_i$ are coincident, which is isomorphic to $M$ as a variety.

\begin{example}
As a variety, $(M / D)^1$ is isomorphic to $M$ and $(M / D)^2$ is isomorphic to the blow-up $\Bl_{D \times D} (M \times M)$. The small diagonal $\Delta_\rel \subseteq (M / D)^2$ is the proper transform of the standard diagonal. In general, we have the natural small diagonal morphism
\begin{align*}
    \iota_{\Delta_\rel} \colon M \cong (M / D)^1 \xrightarrow{\sim} \Delta_\rel \subseteq (M / D)^s. 
\end{align*}
\end{example}

For any subset $S \subseteq \{1, \cdots, r\}$ of cardinality $s$, the moduli space $\oM'_{g, r} (M / D, \beta, \vect{\eta})$ admits a canonical evaluation
\begin{align*}
    \ev_S \colon \oM'_{g, r} (M / D, \beta, \vect{\eta}) \to (M / D)^s,
\end{align*}
which is well-defined by the definition of a relative stable map (the markings are never mapped to the relative divisor).
Suppose $\halpha$ is a partition of length $\hell$. For $\gamma \in H^\ast (M, \bQ)$, let 
\begin{align*}
    \gamma \cdot \Delta_\rel \coloneqq \iota_{\Delta_\rel \ast} (\gamma) \in H^\ast \left((M / D)^{\hell}, \bQ\right).
\end{align*}
We define the relative descendent insertion $\tau_{[\halpha]} (\gamma)$ by
\begin{align}\label{tau_part_rel}
    \tau_{[\halpha]} (\gamma) \coloneqq \psi_1^{\halpha_1 - 1} \cdots \psi_{\hell}^{\halpha_{\hell} - 1} \cdot \ev_{\{1, \cdots, \hell\}}^\ast (\gamma \cdot \Delta_\rel)
\end{align}
Let $\Omega_M (\log D)$ denote the locally free sheaf of differentials with logarithmic poles along $D$. The logarithmic tangent bundle $T_M (- \log D)$ is the dual of $\Omega_M (\log D)$.
For the relative geometry $M / D$, the elements of $\tsK$ also act on $H^\ast (M, \bQ)$ via the substitution $c_i = c_i (T_M (- \log D))$ instead of the substitution $c_i = c_i (T_M)$ used in the absolute case. Then, for even cohomology classes $\gamma_j \in H^{2 \ast} (M, \bQ)$, we define 
\begin{align*}
    \overline{\tau_{\alpha_1 - 1} (\gamma_1) \cdots \tau_{\alpha_{\ell} - 1} (\gamma_{\ell})} = \sum_{\substack{P \text{ set partitions} \\ \text{of } \{1, \cdots, \ell \}}} \prod_{S \in P}  \sum_{0 < |\halpha| \leqslant |\alpha_S|} \tau_{[\halpha]} \left(\tsK_{\alpha_S, \halpha} \cdot \prod_{j \in S} \gamma_j \right)
\end{align*}
as before via \eqref{tau_part_rel} instead of \eqref{tau_part}. In the presence of odd cohomology classes, a sign must be included which is similar to the absolute case (see Remark \ref{sign_rmk}).

Now, we can state the conjectural relative descendent GW/PT correspondence \cite{MNOPII, PP14, PP17}.

\begin{conjecture}\label{relGWPTconj}
Suppose that $D$ is connected and $\alpha=(\alpha_1,\dots,\alpha_r)$ a partition.
For (even) classes $\gamma_j \in H^\ast (M, \bQ)$, $1\leqslant j\leqslant r$, we have
\begin{align*}
    &(-q)^{- \sC_\beta^M / 2} \ZPT \left(M / D; q \tmid \tp_{\alpha_1 - 1} (\gamma_1) \cdots \tp_{\alpha_r - 1} (\gamma_r) \tmid \eta\right)_\beta \\
    ={} &(- i u)^{\sC_\beta^M + \ell (\eta) - |\eta|} \ZGWp \left(M / D; u \tmid \overline{\tau_{\alpha_1 - 1} (\gamma_1) \cdots \tau_{\alpha_r - 1} (\gamma_r)} \tmid \eta \right)_\beta
\end{align*}
under the variable change $e^{i u} = -q$.
\end{conjecture}
The variable change is well-defined assuming Conjecture \ref{relrat_conj}.

\subsection{The degeneration formulas}
Let $W$ be a smooth variety of dimension four and $B$ a smooth irreducible curve with a distinguished point $\mathbf{o} \in B$. 

\begin{definition}
A flat projective morphism $\pi \colon W \to B$ is a {\em simple degeneration} if the following conditions are satisfied:
\begin{enumerate}[(a)]
    \item The morphism $\pi$ has smooth fibers over $ B \setminus \{\mathbf{o}\}$;
    \item The special fiber is the union
    \[
        W_{\mathbf{o}} = M_0 \cup M_1 \cup \cdots \cup M_k
    \] 
     of smooth irreducible components such that for each $1 \leqslant i \leqslant k$, the nonempty intersection $D_i \coloneqq M_0 \cap M_i$ is a smooth connected divisor. Moreover, $M_1, \cdots, M_k$ are pairwise disjoint.
\end{enumerate}
\end{definition}
This definition is a special case of \cite[Def.\ 1.1]{LW15}.

Let $M$ denote a fixed general fiber $W_b$ and \[D \coloneqq \sum_i D_i.\] We will also denote $\pi$ briefly as $M \leadsto M_0 \cup_D (M_1 \sqcup \cdots \sqcup M_k)$.
We write
\begin{align*}
    \iota \colon M \to W, \quad \iota_0 \colon M_0 \to W, \quad \iota_1 \colon M_1 \sqcup \cdots \sqcup M_k \to W
\end{align*}
for inclusions. 

The degeneration formulas express the absolute invariants of $M$ via the relative invariants of $(M_0, D)$ and $(M_1 \sqcup \cdots \sqcup M_k, D)$. We state for completeness the formulas in both GW and PT theories.

\begin{theorem}\label{thm:deg-gw}
Suppose $M$ is a smooth projective threefold. Suppose $\gamma_1,\cdots,\gamma_r$ are even cohomology classes on the total space $W$. 
For a nonzero class $\beta^\prime \in \NE (W)$, we have
\begin{align*}
    &\sum_{\substack{
    \beta \in \NE (M) \\ 
    \iota_\ast \beta = \beta^\prime}} \ZGWp \left(M ;u\tmid \overline{\tau_{\alpha_1-1}(\gamma_1)\cdots\tau_{\alpha_r-1}(\gamma_r)} \right)_{\beta} \\
    ={}&
    \sum \ZGWp \left(M_0 / D; u \tmid \overline{\prod_{i\in I_0}\tau_{\alpha_i-1}(\gamma_i)} \tmid \eta_1, \cdots, \eta_k \right)_{\beta_0}\cdot \\
    &\quad\quad\prod_{j = 1}^k 
    \fz(\eta_j) u^{2 \ell(\eta_j)} \ZGWp \left(M_j/ D_j;u \tmid \overline{\prod_{i\in I_j}\tau_{\alpha_i-1}(\gamma_i)} \tmid \eta_j^{\vee} \right)_{\beta_j}
\end{align*}
where the summation on the second line runs over
\begin{enumerate}[(a)]
    \item splittings 
    \begin{align}\label{eq:curve-splitting-general}
        \iota_{0 \ast} \beta_0 + \iota_{1 \ast} (\sum_{i=1}^k\beta_i) = \beta^\prime=\iota_*\beta
    \end{align} such that $\beta_0 \cdot D_i = \beta_i \cdot D_i$,
    \item partitions $I_0\coprod \cdots \coprod I_k=\{1,2, \cdots, r\}$, and 
    \item cohomology weighted partitions $\eta_1,\cdots,\eta_k$ such that $|\eta_i|=\beta_i \cdot D_i$ with respect to a fixed basis of $H^\ast (D_i, \bQ)$ for $1 \leqslant i \leqslant k$.
\end{enumerate} 
\end{theorem}
See for example \cite{Li02} and \cite[p.403]{PP17}. Similarly, the formula without bars, namely without applying the universal transformation to descendent insertions, also holds.

For the degeneration formulas in symplectic geometry, see \cite{li2001symplectic,IP04,FZ20}.

\begin{theorem}\label{thm:deg-pt}
With notation as in Theorem~\ref{thm:deg-gw}, we have
\begin{align*}
    &\sum_{\substack{\beta \in \NE (M) \\ \iota_\ast \beta = \beta^\prime}} \rZ_\PT \left(M; q \tmid 
    {\tau_{\alpha_1-1}(\gamma_1)\cdots\tau_{\alpha_r-1}(\gamma_r)} \right)_{\beta '} \\
    ={}&
    \sum \rZ_\PT \left(M_0 / D; q \tmid 
    {\prod_{i\in I_0}\tau_{\alpha_i-1}(\gamma_i)} \tmid \eta_1, \cdots, \eta_k \right)_{\beta_0} \cdot\\
    &\quad \quad\prod_{j = 1}^k (- 1)^{|\eta_j| - \ell (\eta_j)} \fz(\eta_j) q^{- |\eta_j|}
    \rZ_\PT \left(M_j/ D_j; q \tmid 
    {\prod_{i\in I_j}\tau_{\alpha_i-1}(\gamma_i)} \tmid \eta_j^{\vee} \right)_{\beta_j}
\end{align*}
where the summation on the second line runs over the same index set in Theorem \ref{thm:deg-gw}. 
\end{theorem}
See for example \cite{LW15} and \cite[p.2761]{PP14}. For the proof of a version of the statement, see \cite[Thm.\ 6.12]{Lin23}. For the parallel statement in DT-theory, see \cite{MNOPII}. Many cases have been proven in \cite[Thm.\ 1.4]{LW15} and \cite[Thm.\ 1.2]{Zhou18}.

Given a splitting \eqref{eq:curve-splitting-general}, we have the following constraints by adjunction for $M \leadsto M_0 \cup_D (M_1 \sqcup \cdots \sqcup M_k)$, which are similar to \cite[Lem.\ 2.2]{HL12}:
\begin{align}\label{eq:dim_constr}
    \sC_\beta^M = \sC_{\beta_0}^{M_0} + \sum_{i = 1}^k (\sC_{\beta_i}^{M_i} - 2 \beta_i \cdot D_i), \\
    \label{eq:int-sing-div}\beta_0 \cdot D_i = \beta_i \cdot D_i \quad\mbox{for } 1 \leqslant i \leqslant k.
\end{align}
This will be important for our arguments.

For notational convenience, we set 
\begin{align*}
    \rZ_\GW^\prime (M_i / D_i; u)_0=\rZ_\PT (M_i / D_i; q)_0 = 1
\end{align*} 
for each $1 \leqslant i \leqslant k$ as a convention. This will appear when the curve misses an irreducible component $M_i$ of the degeneration in the application of degeneration formulas. 

We conclude this section with a well-known lemma. For the convenience of the reader, we provide an argument here.

\begin{lemma} \label{mono_inj}
If the monodromy on $H_2 (M, \bQ)$ around $\mathbf{o} \in B$ is trivial, then $\iota_\ast \colon H_2 (M, \bQ) \to H_2 (W, \bQ)$ is injective and so is $\iota_\ast \colon \NE (M) \to \NE (W)$.
\end{lemma}

\begin{proof}

By hypothesis and the local invariant cycle theorem (see \cite{BBD, Clemens77} or \cite[Thm.\ 4.18]{VoisinII}), the restriction map $\iota^\ast \colon H^2 (W, \bQ) \to H^2 (M, \bQ)$ is surjective. According to the universal coefficient theorem, it follows that 
\[
    \iota_\ast = (\iota^\ast)^\vee \colon H_2 (M, \bQ) \cong H^2 (M, \bQ)^\vee \to H^2 (W, \bQ)^\vee \cong H_2 (W, \bQ)
\]
is injective. Note that $\iota_\ast$ preserves effective cycles by the definition of the pushforward.
\end{proof}

\section{Main Theorem}\label{sec:proof}
We first review the definition of conifold transitions. Let $\pi\colon \fX \to \Delta$ be a projective flat morphism from a smooth fourfold $\fX$ to the unit disk $\Delta$ in $\bC$ and $X$ be a general fiber of it.
\begin{defn}\label{defn:con-trans}
    Suppose the central fiber $\Bar{X}=\fX_0$ of $\pi$ has ordinary double points $\{p_1,\cdots,p_k\}$ as singularities. The morphism $\pi$ together with a projective small resolution $\psi\colon Y\to \Bar{X}$ is a {\em (projective) conifold transition}. We denote it as $X \nearrow Y$. 
\end{defn}

Let $X\nearrow Y$ be a conifold transition and use the notation in the definition.
The following is the main result of the paper.

\begin{thm}\label{thm:main}
    Suppose that $\beta\in \NE(X)$ is a nonzero class and $\alpha=(\alpha_1,\dots, \alpha_r)$ a fixed partition. Asuume $\gamma_i\in H^{\ast}(\fX)$, $i=1,\dots ,r$, are fixed even cohomology classes and if $\gamma_i\in H^0(\mathfrak{X})$, then $\alpha_i=1$. 
    \begin{enumerate}[(a)]
        \item\label{thm:mainA} If Conjecture~\ref{rat_conj} holds for Y, then it holds for $X$ and descendent insertions 
        \begin{align}\label{eq:ins-total}
            \gamma_{i|X},\quad i=1,\dots ,r.
        \end{align}
        \item\label{thm:mainB} If furthermore the GW/PT correspondence, namely Conjecture~\ref{GWPTconj}, holds for $Y$, then it holds for $X$ with descendent insertions \eqref{eq:ins-total}, i.e.,
        \begin{align*}
        &(- q)^{- \sC_\beta / 2} \ZPT \left(X; q \tmid \tp_{\alpha_1 - 1} (\gamma_{1|X})  \cdots \tp_{\alpha_r - 1} (\gamma_{r|X}) \right)_\beta \\
        ={}&(- iu)^{\sC_\beta} \ZGWp \left(X; u\tmid \overline{\tau_{\alpha_1 - 1} (\gamma_{1|X}) \cdots \tau_{\alpha_r - 1} (\gamma_{r|X})}  \right)_\beta 
        \end{align*}
    \end{enumerate}
\end{thm}
The strategy of the proof is to put the smoothing $X$ and the resolution $Y$ into two different degenerations and apply the degeneration formulas.

\begin{remark}\label{rmk_inv}
Since $\oX$ has only ordinary double points, the monodromy acts trivially on $H^q (X, \bZ)$ for $q < 3$ (cf.\ \cite[Cor.\ 2.17]{VoisinII}), and thus for $q > 3$ by Poincar\'e duality. By the local invariant cycle theorem (see \cite{BBD, Clemens77} or \cite[Thm.\ 4.18]{VoisinII}), the restriction map $H^q (\fX) \to H^q (X)$ is surjective for $q \neq 3$.
\end{remark}

\subsection{Two simple degenerations}\label{sim_deg_subsec}
Let $\tY=\Bl_{p_1, \cdots, p_k}\oX$.
We have the following diagram
\begin{equation*}
    \begin{tikzcd}
        & Y \arrow[d, "\psi"] & \tY\arrow[l,"\phi",swap] \arrow[ld
        ]\\
        X & \oX. \arrow[l,"{\rm sm.}"',rightsquigarrow] & 
    \end{tikzcd}
\end{equation*}
We will consider two simple degenerations $\cX \to \Delta$ and $\cY \to \Delta$. 
Special fibers of both $\cX$ and $\cY$ contain the blow-up $\tY$.

For the degeneration $\fX \to \Delta$, there exists a semi-stable degeneration\footnote{The semi-stable reduction $\cX$ does not require the existence of a small resolution of $\fX_0 = \oX$.}$\cX \to \Delta$ via a degree two base change and blow-ups of $\fX$ (for the construction, see \cite{LY06}). 
The special fiber is \[\cX_0=\tY \cup Q_1\cup  \cdots\cup Q_k\] where $Q_i$ is a smooth quadric threefold in $\bP^4$. Let $E_i$ be the exceptional divisors of  the blow-up $\tY$. Then $E_i \cong \bP^1 \times \bP^1$ is a hyperplane section of $Q_i$ in $\bP^4$. 
The blow-up $\tY$ intersects $Q_i$ transversally along $E_i$ and $Q_i \cap Q_j = \varnothing$ for all $i \neq j$. Set $E = \sum_i E_i$. We also denote the degeneration $\cX \to \Delta$ as \begin{align}\label{eq:ss-deg}
    X \leadsto \tY \cup_E (\sqcup_i Q_i).
\end{align}
According to the adjunction formula, $K_{Q_i} = (K_{\bP^4} + Q_i)|_{Q_i} = \cO_{Q_i} (- 3)$ and 
\begin{align}\label{eq:x-local-c1}
    c_1 (T_{Q_i}) = 3 [E_i].
\end{align}

The other degeneration is the deformation to the normal cone $\cY \to \Delta$ which is the composition of the blow-up 
\begin{align*}
    \Phi \colon \cY \coloneqq \Bl_{\sqcup_i C_i \times \{0\}} (Y \times \Delta) \to Y \times \Delta
\end{align*}
with the projection to $\Delta$. Here, each $C_i \coloneqq \psi^{- 1} (p_i)$ is the exceptional $(-1, -1)$-curve of $\psi$. 
The special fiber is \[\cY_0=\tY\cup \tE_1\cup  \cdots\cup \tE_k\] where each $\tE_i$ is isomorphic to $P (\cO_{\bP^1}(- 1)^{\oplus 2} \oplus \cO_{\bP^1})$. Note that $\tY$ is also the blow-up of $Y$ along the exceptional curves $C_i$'s. The transverse intersection $E_i = \tY \cap \tE_i$ is now regarded as the infinity divisor of $\pi_i \colon \tE_i \to C_i \cong \bP^1$.
We also denote the degeneration as
\begin{align}\label{eq:def-norm}
    Y \leadsto \tY \cup_E (\sqcup_i \tE_i).
\end{align}

We include the following well-known fact about the associated monodromies of the degenerations $\cX / \Delta$ and $\cY / \Delta$.

\begin{lemma}\label{mono_triv}

The monodromy of $\cX /\Delta$ (resp.\ $\cY/\Delta$) around the special fiber $\cX_0$ (resp.\ $\cY_0$) acts trivially on $H_2 (X, \bZ)$ (resp.\ $H_2 (Y, \bZ)$). 
\end{lemma}

\begin{proof}
Since the monodromy on $H_2 (\fX_t, \bZ)$ ($t \neq 0$) around $0 \in \Delta$ is trivial (see \cite[Lem.\ 1.4]{LY06} or \cite[Cor.\ 2.17]{VoisinII}), so is the monodromy on $H_2 (\cX_t, \bZ)$ around $0 \in \Delta$. The same holds for the family $\cY \to \Delta$ since the punctured family is trivial.
\end{proof}

Because $H_2 (X, \bZ)$ has trivial monodromy over $\Delta \setminus \{ 0 \}$,
it is canonically isomorphic to $ H_2 (\oX, \bZ)$ (see \cite[(1.28)]{Clemens83} or \cite[Lem.\ 1.4]{LY06}). Moreover, we have the following exact sequence (see \cite[(1.29)]{Clemens83})
\begin{align*}
    \bigoplus_{i = 1}^k \bZ [C_i] \to H_2 (Y, \bZ) \xrightarrow{\psi_\ast} H_2 (\oX, \bZ) \cong H_2 (X, \bZ) \to 0.
\end{align*}
We set $\phi = \Phi |_{\tY}\colon \tY \to Y$, which is regarded as  the blow-up $\tY \to Y$ of $Y$ along $C_i$'s (via the identification $Y \cong Y \times \{0\} \subseteq Y \times \Delta$).
Note that $\phi$ induces an injective Gysin homomorphism (cf.\ \cite[Rmk.\ 7.29]{VoisinI})
\begin{align}\label{gysin}
    \phi^{!} \coloneqq \mathrm{PD}_{\tY} \circ \phi^\ast \circ \mathrm{PD}_Y \colon H_2 (Y, \bZ) \rightarrowtail H_2 (\tY, \bZ)
\end{align}
where $\mathrm{PD}_{(-)}$ is the Poincar\'e duality isomorphism. 
Moreover,
\begin{align}\label{gysin_im}
    \Img (\phi^{!}) = \{ \Tilde{\beta} \in H_2 (\tY, \bZ) \mid \Tilde{\beta} \cdot E_j = 0 \text{ for }  1 \leqslant j \leqslant k \}.
\end{align}
Under the identification $H_2 (\oX, \bZ) \cong H_2 (X, \bZ)$, we have the following commutative diagram
\begin{equation}\label{eq:diag-h2}
    \begin{tikzcd}[column sep=tiny]
        H_2 (Y, \bZ) \ar[d, "\psi_\ast", two heads] \ar[rrr, "\phi^!", tail] &&& H_2 (\tY, \bZ) \ar[d, "\iota_{0 \ast}"] \\
        H_2(\oX, \bZ) \ar[r, phantom, "\isom"] & H_2(X,\bZ) \ar[rr, "\iota_\ast", tail] && H_2 (\cX, \bZ).
    \end{tikzcd}    
\end{equation}

\subsection{GW/PT correspondence}

We prove the following key proposition for both GW and PT-invariants, which relates invariants over $X$ and $Y$ to those over $\Tilde{Y}$.
This is done by applying the degeneration formulas to the degenerations \eqref{eq:ss-deg} and \eqref{eq:def-norm}.

\begin{proposition}\label{prop:key}
Let $X \nearrow Y$ be a conifold transition of smooth projective threefolds and keep the notation as in Section \ref{sim_deg_subsec}. Suppose that $\alpha=(\alpha_1,\dots, \alpha_r)$ is a fixed partition.
\begin{enumerate}[(a)]
    \item\label{Xabrel} Let $\beta \in \NE (X)$ be a nonzero class and $\gamma_1,\cdots,\gamma_r\in H^{\ast}(\fX)$ be even classes. Suppose if $\gamma_i\in H^0(\mathfrak{X})$, then $\alpha_i=1$. 
    Then we have\footnote{In the four displayed equalities in this proposition, the curve classes do not intersect the corresponding divisors. The weighted cohomology partitions are empty. So, we have omitted the weighted partitions from the relative generating functions, as indicated in \S \ref{rel_theories_subsusec}.}
          \begin{align*}
              \rZ^\prime_\GW \left(X; u\tmid \overline{\prod_{i=1}^r\tau_{\alpha_i-1}(\gamma_i)} \right)_{\beta} 
              &= \sum_{\substack{
              \psi_\ast \beta_Y = \beta}} \rZ^\prime_\GW \left(\tY / E; u \tmid \overline{\prod_{i=1}^r\tau_{\alpha_i-1}(\gamma_i)}\right)_{\phi^{!} \beta_Y},\\
              \rZ_\PT \left(X;q \tmid \prod_{i=1}^r\tau_{\alpha_i-1}(\gamma_i) \right)_{\beta} 
              &= \sum_{\substack{
              \psi_\ast \beta_Y = \beta}} \rZ_\PT \left(\tY / E; q \tmid \prod_{i=1}^r\tau_{\alpha_i-1}(\gamma_i)\right)_{\phi^{!} \beta_Y}
          \end{align*}
          where the summations are finite. 
          Here, by abuse of notation, $\gamma_i$ on the left is viewed as a class on $X$ as a pull-back via the inclusion $X\hookrightarrow \mathfrak{X}$, and $\gamma_i$ on the right is viewed as a class on $\Tilde{Y}$ via the obvious maps $\Tilde{Y}\hookrightarrow \mathcal{X}_0 \to \mathfrak{X}$. 
    \item \label{Yabrel} 
          Let $\tC_{j}$ be the rational curve in 
          $\tE_j$ via the inclusion $\cO_{C_j} \to N_{C_j / Y} \oplus \cO_{C_j}$. Let $\beta_Y \in \NE (Y)$ be a nonzero class and $\gamma_1,\cdots,\gamma_r\in H^{*}(Y)$ be even classes. Suppose either their restrictions to the exceptional curves $C_j$, $j=1,\dots, k$, are zero, or if $\gamma_i\in H^0(\mathfrak{X})$, then $\alpha_i=1$. 
          Then we have
          \begin{align*}
              &\rZ_\GW^\prime \left(Y;u \tmid \overline{ \prod_{i=1}^r\tau_{\alpha_i-1}(\gamma_i)}\right)_{\beta_Y} \\
              =&\sum
              \rZ_\GW^\prime \left(\tY / E; u \tmid \overline{ \prod_{i=1}^r \tau_{\alpha_i-1}(\gamma_i) }\right)_{\phi^{!} {\beta_Y^\prime}} 
              \prod_{j = 1}^k 
              \rZ_\GW^\prime \left(\tE_j / E_j; u \right)_{m_j [\tC_j]}\mbox{ and}
              \\
              &\rZ_\PT \left(Y ;q \tmid \prod_{i=1}^r\tau_{\alpha_i-1}(\gamma_i)\right)_{\beta_Y}\\
              =&\sum
               \rZ_\PT \left(\tY / E; q  \tmid \prod_{i=1}^r \tau_{\alpha_i-1}(\gamma_i) \right)_{\phi^{!} {\beta_Y^\prime}} 
               \prod_{j = 1}^k 
               \rZ_\PT \left(\tE_j / E_j; q\right)_{m_j [\tC_j]} .
          \end{align*}
          The summations are over curve splittings in $H_2(\mathcal{Y})$, omitting the obvious push-forwards, 
          \begin{align*}
              \beta_Y = \phi^{!} {\beta_Y^\prime} + \sum_{j=1}^k m_j [{\tC_j}]
          \end{align*}
          where $\beta_Y^\prime \in \NE (Y)$ and $ m_i \in \bZ_{\geqslant 0}$.
\end{enumerate}
Furthermore, the summations are finite.
\end{proposition}
\begin{proof}
\eqref{Xabrel} By the string equation for GW-invariants (\cite[\S 26.3]{HKK+03}), if the descendent insertions involve $\tau_0(1)$, then the calculation can be reduced to one with fewer insertions, removing this $\tau_0(1)$. 
On the PT-invariants side, by the string equation for PT-invariants (\cite[Prop. 1.1]{Mor22}), if the descendent insertions involve $\tau_0(1)$, then the PT-invariant is zero. 
Thus we can assume $\gamma_i \in H^{>0}(\mathfrak{X})$ and only need to consider cohomology classes of degrees $2$, $4$, and $6$.
    We consider the degeneration \eqref{eq:ss-deg}. Since the map from $Q_j$ to $\mathfrak{X}$ factors through $\{p_j\}\hookrightarrow \mathfrak{X}$ \cite[\S 1]{LY06}, the pullback of $\gamma_i\in H^{>0}(\mathfrak{X})$ to $Q_j$ is zero for all $i$ and $j$. Thus the degeneration formulas will simplify so that there are no descendent insertions over $Q_j$.
    
    For GW-invariants, we apply Theorem \ref{thm:deg-gw}. 
    First notice that it is enough to prove the corresponding statement without bars, namely without applying the universal transformation on descendent insertions. 
    On the other hand, according to Lemmas~\ref{mono_inj} and \ref{mono_triv}, we have $\NE(X)\hookrightarrow \NE(\fX)$. Then the degeneration formula is simplified to 
    \begin{align*}
    &
    \rZ_\GW ^\prime\left( X ;u \tmid {\tau_{\alpha_1-1}(\gamma_1)\cdots\tau_{\alpha_r-1}(\gamma_r)} \right)_{\beta} \\
    ={}&
    \sum
    \rZ_\GW^\prime \left(\Tilde{Y} / E; u \tmid {\prod_{i=1}^r\tau_{\alpha_i-1}(\gamma_i)} \tmid \eta_1, \cdots, \eta_k \right)_{\Tilde{\beta}}\cdot\\ 
    &\quad \quad\prod_{j = 1}^k \fz(\eta_j) u^{2 \ell(\eta_j)}
    \rZ_\GW^\prime \left(Q_j/ E_j; u  \tmid \eta_j^{\vee} \right)_{\beta_j}.
\end{align*}
Here, the summation is over the curve splittings\footnote{Here and in \eqref{eq:curve-splitting-ycal}, we have omitted various push-forward symbols to simplify the notation.} in $H^2(\mathcal{X})$ 
\begin{align}\label{eq:curve-splitting-gwpt}
\beta=\Tilde{\beta}+\sum\beta_j\end{align}
and cohomology weighted partitions $\eta_1,\cdots,\eta_k$ such that 
\begin{align}\label{eq:compatible-int-divisor}
    \Tilde{\beta}\cdot E_j=\beta_j\cdot E_j=|\eta_j|.
\end{align}
    Since there are no descendant insertions over $Q_j$, 
    according to the virtual dimensions of the moduli spaces of relative stable maps \eqref{rel_vdGW},
    \begin{align*}
        \sC^{{Q_j}}_{\beta_j}+\ell(\eta_j^\vee)-|\eta_j^\vee| = \sC^{{Q_j}}_{\beta_j}+\ell(\eta_j)-|\eta_j|=0.
    \end{align*}
    Moreover, according to \eqref{eq:x-local-c1}, we have 
    \begin{align}\label{eq:int-div-0}
    \ell(\eta_j)=0=\beta_j\cdot E_j=|\eta_j|\quad\mbox{for }1\leqslant j\leqslant k.
    \end{align}
    In particular, all weighted cohomology partitions are empty.
    Since $E_j\subseteq Q_j$ is a hyperplane section, $\beta_j\cdot E_j=0$ implies that $\beta_j=0$. 
    Then, according to \eqref{gysin_im} and \eqref{eq:diag-h2}, there is a unique $\beta_Y\in H_2(Y,\mathbb{Z})$ such that $\phi^{!}\beta_Y=\Tilde{\beta}$ and $\psi_\ast\beta_Y=\beta$. 
    Thus, the degeneration formula simplifies to the desired form, 
    using the diagram \eqref{eq:diag-h2}.
    
    Similarly, for PT-invariants, we apply Theorem~\ref{thm:deg-pt}, obtaining    
    \begin{align*}
    &
    \rZ_\PT \left(X; q \tmid {\tau_{\alpha_1-1}(\gamma_1)\cdots\tau_{\alpha_r-1}(\gamma_r)} \right)_{\beta} \\
    = &
    \sum
    \rZ_\PT \left(\Tilde{Y} / E ;q \tmid {\prod_{i=1}^r\tau_{\alpha_i-1}(\gamma_i)} \tmid \eta_1, \cdots, \eta_k \right)_{\Tilde{\beta}}\cdot \\
    &\quad\quad\prod_{j = 1}^k (- 1)^{|\eta_j| - \ell (\eta_j)} \fz(\eta_j) q^{- |\eta_j|}
    \rZ_\PT \left(Q_j/ E_j; q \tmid \eta_j^{\vee} \right)_{\beta_j}. 
\end{align*}
The summation is over curve splittings \eqref{eq:curve-splitting-gwpt} and cohomology weighted partitions $\eta_1,\cdots,\eta_k$ satisfying \eqref{eq:compatible-int-divisor}.
According to \eqref{eq:rel-pt-vdim}, 
    \begin{align*}
        \sC_{\beta}^X&=\sum_{i=1}^r(\alpha_i-1+\deg \gamma_i)\quad\mbox{and} \\
        \sC_{\Tilde{\beta}}^{\Tilde{Y}}&=\sum_{i=1}^r(\alpha_i-1+\deg \gamma_i)+ \frac{1}{2}\sum_{j=1}^k\deg_{\mathbb{R}} C_{\eta_j}.
    \end{align*}
    On the other hand, according to \eqref{eq:dim_constr}, we have
    \begin{align*}
        \sC^X_{\beta}=\sC^{\Tilde{Y}}_{\Tilde{\beta}} +\sum_{j=1}^k(\sC^{Q_j}_{\beta_j}-2\beta_j\cdot E_j).
    \end{align*}
    Combined with \eqref{eq:x-local-c1}, we have $\sum_{j=1}^k(\deg_{\mathbb{R}} C_{\eta_j}/2+|\eta_j|)=0$. Since all the summands are non-negative, we again have \eqref{eq:int-div-0} and 
    $\beta_j=0$.
    Then the equality follows from the degeneration formula.
    
    Note that the finiteness of the sums in the right-hand side of \eqref{Xabrel} follows from that $\phi^{!} \beta_Y \notin \NE(\tY)$ for all but finitely many $\beta_Y$, for a proof see \cite[Prop. 2.1, p.511]{LLW18}.
    
    \vspace{6pt}
    \noindent \eqref{Yabrel} Again by the string equations, we only need to consider insertions of positive degrees. 
    We consider the pullback of the descendent insertion $\gamma_i$ to $\mathcal{Y}$ via the composition
    \begin{align*}
        \mathcal{Y}\xrightarrow{\Phi} Y\times \Delta \xrightarrow{\operatorname{pr}} Y. 
    \end{align*}
    Since the induced map $\Tilde{E}_j\to Y$ factors through $C_j\subseteq Y$, the pullback to $\Tilde{E}_j$ is zero.
    
    We first consider GW-invariants. Again, it is enough to prove the corresponding statement without bars.
    Applying Theorem~\ref{thm:deg-gw} and Lemmas \ref{mono_inj} and \ref{mono_triv} to the degeneration \eqref{eq:def-norm},
    we obtain
    \begin{align*}
    &\rZ_\GW ^\prime\left( Y; u \tmid {\tau_{\alpha_1-1}(\gamma_1)\cdots\tau_{\alpha_r-1}(\gamma_r)} \right)_{{\beta_Y}} \\
    = {}&
    \sum
    \rZ_\GW^\prime \left(\Tilde{Y} / E; u \tmid {\prod_{i=1}^r\tau_{\alpha_i-1}(\gamma_i)} \tmid \eta_1, \cdots, \eta_k \right)_{\Tilde{\beta}} 
    \prod_{j = 1}^k \fz(\eta_j) u^{2 \ell(\eta_j)}
    \rZ_\GW^\prime \left(\Tilde{E}_j/ E_j; u  \tmid \eta_j^{\vee} \right)_{\beta_j}.
    \end{align*}
    The summation is over the curve splittings
    \begin{align}\label{eq:curve-splitting-ycal}
        \beta_Y=\Tilde{\beta}+\sum_{j=1}^k\beta_j,
    \end{align}
    in $H_2(\mathcal{Y})$ and cohomology weighted partitions $\eta_1,\cdots,\eta_k$ satisfying 
    \eqref{eq:compatible-int-divisor}. 
    Note that $\sC_{\beta_j}^{\tE_j} = 3 \beta_j \cdot E_j = 3 |\eta_j|$ 
    because of $\tE_j \cong P(\cO_{\bP^1}(- 1)^{\oplus 2} \oplus \cO_{\bP^1})$ and \eqref{loccur_T}.
    Similar to the previous part, we can obtain \eqref{eq:int-div-0}. 
    Because $[E_i] = c_1 (\cO_{\tE_i} (1))$, 
    $\beta_j = m_j [\tC_j]$ for some $m_j \in \bZ_{\geqslant 0}$. 
    Also, because of \eqref{eq:compatible-int-divisor}, there is a unique $\beta_Y^\prime \in \NE (Y)$ such that $\Tilde{\beta} = \phi^! \beta_Y^\prime$ by \eqref{gysin_im}. 
    Hence, \eqref{Yabrel} is proven for GW-invariants.

    For PT-invariants, we apply Theorem~\ref{thm:deg-pt} and Lemmas \ref{mono_inj} and \ref{mono_triv} to the degeneration \eqref{eq:def-norm}, 
    obtaining 
    \begin{align*}
    &
    \rZ_\PT \left(Y; q \tmid {\tau_{\alpha_1-1}(\gamma_1)\cdots\tau_{\alpha_r-1}(\gamma_r)} \right)_{{\beta_Y}} \\
    ={} &
    \sum
    \rZ_\PT \left(\Tilde{Y} / E; q  \tmid {\prod_{i=1}^r\tau_{\alpha_i-1}(\gamma_i)} \tmid \eta_1, \cdots, \eta_k \right)_{\Tilde{\beta}}\cdot\\ 
    &\quad\quad \prod_{j = 1}^k (- 1)^{|\eta_j| - \ell (\eta_j)} \fz(\eta_j) q^{- |\eta_j|}
    \rZ_\PT \left(\Tilde{E}_j/ E_j; q \tmid \eta_j^{\vee} \right)_{\beta_j}. 
    \end{align*}
    The summation is over the curve splittings \eqref{eq:curve-splitting-ycal}
    and cohomology weighted partitions $\eta_1,\cdots,\eta_k$ satisfying 
    \eqref{eq:compatible-int-divisor}.
    We can again deduce \eqref{eq:int-div-0}.
    Thus, $\beta_j$ takes the desired form, similar to the GW case. Then the degeneration formula reduces to the equality in the statement.

    Finally, we note that \eqref{eq:curve-splitting-ycal} is equivalent to the curve splitting $\beta_Y = \phi_\ast \Tilde{\beta} + \sum m_j \pi_{j \ast} [\tC_j]$ in $\NE (Y)$, where $\pi_j \colon \tE_j \to C_j$. By \cite[Cor.\ 1.19 (3)]{KM98}, the summations in the right-hand side of \eqref{Yabrel} are finite.
\end{proof}

Now, we are ready to prove the main theorem.

\begin{proof}[Proof of Theorem~\ref{thm:main}]
Keeping the notation of Proposition \ref{prop:key} and its proof, we first prove the following correspondence without weighted cohomology partitions (cf.\ Conjecture \ref{relGWPTconj})
\begin{align}
    & (-q)^{- \sC_{\Tilde{\beta}}^{\Tilde{Y}} / 2} \ZPT \left(\Tilde{Y}/ E ; q\tmid \tp_{\alpha_1 - 1} (\gamma_{1}) \cdots \tp_{\alpha_r - 1} (\gamma_{r}) \right)_{\Tilde{\beta}} \nonumber\\
    =& (- i u)^{\sC_{\Tilde{\beta}}^{\Tilde{Y}}} \ZGWp \left( \Tilde{Y} / E; u \tmid \overline{\tau_{\alpha_1 - 1} (\gamma_{1}) \cdots \tau_{\alpha_r - 1} (\gamma_{r})} \right)_{\Tilde{\beta}} \label{eq:key-equality}
\end{align}
for all $\Tilde{\beta} \in \operatorname{Im} (\phi^!) \subseteq \NE (\Tilde{Y})$ under the variable change $- q = e^{iu}$. 
Here $\gamma_i$, $i=1,\dots,r$, are pulled back to $\Tilde{Y}$ via the composition $\Tilde{Y} \to Y\to \mathfrak{X}$. Under the map $Y\to \mathfrak{X}$, the curve $C_j$ is mapped to $p_j\in \mathfrak{X}$. Thus, the pullbacks are restricted to zero on $C_j$. 
By Proposition \ref{prop:key} \eqref{Yabrel} and $K_Y \cdot C_j = 0$ for $1 \leqslant j \leqslant k$, we have
\begin{align}
\label{eq:keyformula}
    &\sum_{\beta_Y \in \NE (Y)} (-q)^{- \sC_{\beta_Y} / 2} \rZ_\PT \left(Y ;q \tmid \tau_{\alpha_1 - 1} (\gamma_{1}) \cdots \tau_{\alpha_r - 1} (\gamma_{r})\right)_\beta v^{\beta_Y} \\
    ={}& \left( \sum_{\beta_Y^\prime \in \NE (Y)} (-q)^{- \sC_{\beta_Y'} / 2} \rZ_\PT \left(\tY / E; q \tmid \tau_{\alpha_1 - 1} (\gamma_{1}) \cdots \tau_{\alpha_r - 1} (\gamma_{r})\right)_{\phi^{!} \beta_Y^\prime} v^{\beta_Y^{\prime}} \right) \cdot \nonumber\\
    &\quad\quad \prod_{i = 1}^k \left( \sum_{m_i \geqslant 0} \rZ_\PT (\tE_i / E_i; q)_{m_i [\tC_{i}]} v^{m_i [\Tilde{C_i}]} \right).\nonumber
\end{align}
We have a similar equality for GW-invariants.
By abuse of notation, let 
\begin{align*}
    \PTop \left(Y \tmid \prod_{j=1}^r \tau_{\alpha_j - 1} \left(\gamma_{j}\right) \right) 
    &\coloneqq \sum_{\beta_Y \in \NE (Y)} (-q)^{- \sC_{\beta_Y} / 2} \ZPT \left(Y; q \tmid \prod_{j=1}^r \tau_{\alpha_j - 1} \left(\gamma_{j}\right)\right)_{\beta_Y} v^{\beta_Y}, \\
    \PTexc (Y) &\coloneqq \prod_{i = 1}^k \left( \sum_{m_i \geqslant 0} \ZPT (\tE_i / E_i; q)_{m_i [\tC_{i}]} v^{m_i [\Tilde{C_i}]} \right). 
    \end{align*}
    We define
    \begin{align*}
    \GWop \left(Y \tmid \overline{\prod_{j=1}^r \tau_{\alpha_j - 1} (\gamma_{j})}\right) &\coloneqq \sum_{\beta_Y \in \NE (Y)} (- i u)^{\sC_{\beta_Y}} \ZGWp \left(Y; q \tmid \overline{\prod_{j=1}^r \tau_{\alpha_j - 1} (\gamma_{j})}\right)_{\beta_Y} v^{\beta_Y},
\end{align*}
and $\GWpexc (Y)$ similarly.
Since $\rZ_\PT (\tE_i / E_i; q)_0 =\rZ_\GW^\prime (\tE_i / E_i; u)_0 = 1$, the generating series $\PTexc (Y)$ and $\GWpexc (Y)$ are invertible in $\bQ (\!( q )\!) [\![ \NE (Y) ]\!]$ and $\bQ (\!( u )\!) [\![ \NE (Y) ]\!]$ respectively. Then we may rewrite \eqref{eq:keyformula} as 
\begin{align}\label{PTratio}
    &\frac{\PTop \left(Y \tmid \prod_{j=1}^r \tau_{\alpha_j - 1} \left(\gamma_{j}\right)\right)}{\PTexc (Y)} \\=& \sum_{\beta_Y^\prime \in \NE (Y)} (-q)^{- \sC_{\beta_Y'} / 2} \ZPT \left({\tY} / E; q \tmid \prod_{j=1}^r \tau_{\alpha_j - 1} \left(\gamma_{j}\right)\right)_{\phi^{!} \beta_Y^\prime} v^{ \beta_Y^\prime},\nonumber
\end{align}
and similarly for $\GWop\left(Y\tmid \overline{\prod_{j=1}^r \tau_{\alpha_j - 1} (\gamma_{j})}\right) / \GWpexc (Y)$.

Since $\bR_{\geqslant 0} [\tC_i]$ is an extremal ray of $\NE (\tE_i)$ and $\tC_i$ does not intersect $E_i$, we have 
\begin{align}
    \rZ_\PT (\tE_i / E_i; q)_{m_i [\tC_{i}]} &= \rZ_\PT (\tE_i; q)_{m_i [\tC_{i}]}\quad \mbox{and}\label{eq:etilde-pt}\\
    \rZ_\GW^\prime (\tE_i / E_i; u)_{m_i [\tC_{i}]} &= \rZ_\GW^\prime (\tE_i; u)_{m_i [\tC_{i}]}\label{eq:etilde-gw}
\end{align}
by applying the degeneration formula to the deformation to the normal cone $\tE_i \leadsto \Bl_{\tC_i} \tE_i \cup_{E_i} \tE_i$. 
The rationality result, namely Conjecture \ref{relrat_conj}, holds for \eqref{eq:etilde-pt}
, since $\widetilde{E}_i$ is toric. 
According to Theorem \ref{GWPTloc}, it follows that under the variable change $q = - e^{i u}$ we have $\ZPT (\tE_i / E_i; q)_{m_i [\tC_{i}]} = \ZGWp (\tE_i / E_i; u)_{m_i [\tC_{i}]}$ and thus $\PTexc (Y) = \GWpexc (Y)$.

If Conjecture~\ref{rat_conj} holds for $Y$, then $\ZPT \left({\tY} / E ; q \tmid \prod_{j=1}^r \tau_{\alpha_j - 1} \left(\gamma_{j}\right)\right)_{\phi^{!} \beta_Y}$ is rational. Thus, $\ZPT \left(X ; q \tmid \prod_{j=1}^r \tau_{\alpha_j - 1} \left(\gamma_{j}\right)\right)_\beta$ is rational, according to Proposition~\ref{prop:key} \eqref{Xabrel}. We have proven \eqref{thm:mainA}.  

If both Conjectures~\ref{rat_conj} and \ref{GWPTconj} hold for $Y$, then
\begin{align*}
    \frac{\PTop \left(Y \tmid \prod_{j=1}^r \tau_{\alpha_j - 1} \left(\gamma_{j}\right) \right)}{\PTexc (Y)} = \frac{\GWop\left(Y \tmid \overline{\prod_{j=1}^r \tau_{\alpha_j - 1} \left(\gamma_{j}\right)}\right)}{\GWpexc (Y)},
\end{align*}
and \eqref{eq:key-equality} follows by extracting the coefficient of $v^{ \Tilde{\beta}}$ from both sides of the above equality.

For each nonzero $\beta \in \NE (X)$, the descendent $\GW / \PT$ correspondence for $X$ now follows from Proposition \ref{prop:key} \eqref{Xabrel} and taking the sum of both sides of \eqref{eq:key-equality} over all $\beta_Y \in \NE (Y)$ with $\psi_\ast \beta_Y = \beta$. 
\end{proof}

\section{Applications}\label{sec:app}
We apply our main result to Fano threefolds and double solids. We also raise a question (Question \ref{qst}) about general small transitions of Calabi--Yau threefolds.

\subsection{Fano threefolds via small toric degenerations} 

We fix notation for Fano threefolds, namely smooth projective threefolds with an ample anticanonical line bundle, as follows:
\begin{itemize}

    \item Let $Q_2$ denote a smooth quadric hypersurface in $\bP^4$. Let $B_k$ (resp.\ $V_k$) denote the Fano threefold with Picard number $1$, Fano index $2$ (resp.\ $1$) and anti-canonical degree $(- K)^3 = 8k$ (resp.\ $ k$). 

    \item Let $V_{\rho, n}$ denote the $n$-th entry in the Mori--Mukai list \cite[Table 2-4]{MM82} of Fano threefolds of Picard number $\rho$.
\end{itemize}
Deformation families of Fano threefolds have been completely classified, see \cite{Iskovskih77, Iskovskih78, MM82, MM82err}. In Galkin's thesis (see also \cite{galkin2018small}), he described all conifold transitions from such Fano threefolds to toric weak Fano threefolds.

\begin{theorem}[\cite{galkin2018small}]\label{thm:Galkin}
There are $44$ families of non-toric Fano threefolds $X$ which admits conifold transitions $X \nearrow Y$ to toric threefolds $Y$:
\begin{enumerate}[(a)]
    \item For $\rho (X) = 1$, there are $4$ families: $Q_2, B_4, B_5, V_{22}$.
    \item For $\rho (X) = 2$, there are $16$ families: $V_{2, n}$ where $n = 12, 17$ or $19 \leqslant n \leqslant 32$.
    \item For $\rho (X) = 3$, there are $16$ families: $V_{3, n}$ where $n = 7$ or $10 \leqslant n \leqslant 24$.
    \item For $\rho (X) = 4$, there are $8$ families: $V_{4, n}$ where $1 \leqslant n \leqslant 8$.
\end{enumerate}
\end{theorem}

Applying Theorems \ref{GWPTtoric}, \ref{thm:main}, and \ref{thm:Galkin}, we conclude: 

\begin{cor}\label{cor:non-toric}
Let $X$ be one of the Fano threefolds in Theorem \ref{thm:Galkin}. Then the GW/PT correspondence (Conjecture~\ref{GWPTconj}) holds for $X$ with descendent insertions \eqref{eq:ins-total}.
\end{cor}

\subsection{Double solids}

Double covers of $\bP^3$ with at worst ordinary double point singularities, which obtained the name double solids, were studied by Clemens \cite{Clemens83}. The construction of Clemens has straightforward generalizations to more general Fano threefolds (cf.\ \cite[\S 4.1]{GH88}). 
We still call the resulting double cover a double solid.

To apply Theorem \ref{thm:main}, we need the following proposition, which is probably well-known. For lack of a suitable reference we will give a sketch of a proof.

\begin{prop}\label{dcover}
Suppose $(Z,\cL)$ is one of the following pairs:
\begin{enumerate}[(a)]
    \item $(Z,\cL) = (\bP^3,\cO_{\bP^3} (a))$ for $a = 2, 3 ,4$;
    \item $(Z,\cL) = (\bP^1 \times \bP^2, \cO_{\bP^1}(1) \boxtimes \cO_{\bP^2}(b))$ for $b = 1, 2$;
    \item $(Z,\cL) = (\bP^1 \times \bP^1 \times \bP^1,\cO_{\bP^1} (1)^{\boxtimes 3})$.
\end{enumerate}
Let $Y$ be the zero locus in $Z \times \bP^1$ defined by a general section $s \in H^0 (Z \times \bP^1, \cL \boxtimes \cO_{\bP^1}(2))$, and $X$ the double cover of $Z$ branched along a smooth surface defined by a general section of $H^0 (Z, \cL^{\otimes 2})$. Then there is a conifold transition $X \nearrow Y$ from $X$ to $Y$.
\end{prop}

\begin{proof}[Sketch of proof]
First, the $Y$ is a smooth hypersurface in a product $Z \times \bP^1$ of projective spaces by Bertini’s theorem. Let $x_0$ and $ x_1$ be homogeneous coordinates on $\bP^1$. For the general section $s$, there are sections $s_{ij} \in H^0 (Z, \cL)$ such that 
\begin{align}\label{qua}
    s = \sum_{0 \leqslant i \leqslant j \leqslant 1} s_{ij} x_i x_j.
\end{align}
Let $Y \to \oX \to Z$ be the Stein factorization of the restriction of the projection $ Z \times \bP^1 \to Z$ to $Y$. 
Then $\oX$ is a double cover of $Z$ branched along a surface $B$ defined by the discriminant of the qudratic
equation \eqref{qua} in $x_0$ and $x_1$, given by $s_{01}^2 - 4 s_{00} s_{11} \in H^0 (Z, \cL^{\otimes 2})$. A local computation shows that $B$ is a nodal surface and thus $\oX$ has only ordinary double points. By perturbing the general section of $\cL^{\otimes 2}$ to the discriminant of \eqref{qua}, we get a projective degeneration of $X$ to the double solid $\oX$ and hence get a conifold transition $X \nearrow Y$.
\end{proof}

\begin{remark}
For $\cL = \cO_{\bP^3} (4)$, the $X$ is a Calabi--Yau threefold, which was studied in  \cite[Prop.\ 3]{GH88} (see also Example 5.8 in \cite{Wang18}). Example 1.7 in \cite{CP10} considered the Fano threefold $X$ associated to $\cL = \cO_{\bP^3} (3)$. Table \ref{tab:Fsolid} gives the corresponding Fano threefolds $X$ in Proposition \ref{dcover}. 
\end{remark}

\begin{table}[H]
  \centering
  \begin{tabular}{c|ccccc}
    \toprule
    $\cL$ & $\cO_{\bP^3} (3)$ & $\cO_{\bP^3} (2)$ & $\cO_{\bP^1}(1) \boxtimes \cO_{\bP^2}(2)$ & $\cO_{\bP^1}(1) \boxtimes \cO_{\bP^2}(1)$ & $\cO_{\bP^1} (1)^{\boxtimes 3}$  \\
    \midrule
    $X$ & $V_2$ & $B_2$ & $V_{2,2}$ & $V_{2,18}$ & $V_{3,1}$ \\
    \bottomrule
  \end{tabular}
  \caption{Fano threefolds realized as double covers}
  \label{tab:Fsolid}
\end{table}

\begin{cor}\label{cor:dcover}
Let $X$ be one of the smooth double covers in Proposition \ref{dcover}. Then the GW/PT correspondence (Conjecture~\ref{GWPTconj}) holds for $X$ with descendent insertions \eqref{eq:ins-total}.    
\end{cor}

\begin{proof}[Sketch of proof]
Let $X \nearrow Y$ be the conifold transition as in Proposition \ref{dcover}. If $Y$ is Calabi--Yau, i.e., $\cL = \cO_{\bP^3} (4)$, then $Y$ satisfies the $\GW/\PT$ correspondence by \cite[Thm.\ 1]{PP17} since it is a smooth hypersurface in a product of projective spaces. For the others, we can also degenerate the weak Fano complete intersection threefolds $Y$ by a similar factoring argument in the proof of \cite[Thm.\ 1]{PP17}. Hence the $\GW/\PT$ correspondence holds for $Y$. Then Corollary \ref{cor:dcover} follows from Theorem \ref{thm:main}.
\end{proof}

\begin{remark}\label{rmk:allFano3}
One can apply Theorem \ref{thm:main} inductively to a sequence of conifold transitions $Y_{i - 1} \nearrow Y_i$ $(1 \leqslant i \leqslant r)$ with $Y_0 = X$. In Corollay \ref{cor:non-toric} and \ref{cor:dcover}, we only use \emph{one} conifold transition. For other Fano threefolds $X$, one may find such a sequence with $Y_r$ being toric ($r \geqslant 2$). For example, using Corti--Hacking--Petracci's construction \cite[\S 6.4.1]{CCdS22}, one has 
conifold transitions
\[
    Y_0 \nearrow Y_1 \nearrow Y_2
\]
such that $Y_0$ is the complete intersection $V_8$ of three quadrics in $\bP^6$ and $Y_2$ is toric, even though the correspondence is known to hold for both $Y_0$ and $Y_2$ by the works of Pandharipande and Pixton.
\end{remark}


\subsection{Concluding remark} 
We make a comment on the ratio \eqref{PTratio} for a Calabi--Yau threefold $Y$. If $\psi \colon Y \to \oX$ is a flopping contraction, then the transformation formula of the ratio $\PTop (Y) / \PTexc (Y)$ under flops was proved by Toda (\cite[Thm.\ 1.2]{Toda10} \& \cite[Thm.\ 1.3]{Toda13}) and Calabrese \cite{Calabrese16}, where the generating series of exceptional curves of $\psi$ is defined by
\begin{align*}
    \PTexc (Y) = \sum_{\substack{\beta_Y \in \NE (Y) \\ \psi_\ast \beta_Y = 0}} \ZPT (Y; q)_{\beta_Y} v^{\beta_Y}.
\end{align*}
Assume $\oX$ is smoothable, and let $X$ denote a smoothing of it. Such an extremal transition $X \nearrow Y$ is called \emph{small} in \cite[Def.\ 6.1]{Wang18}.  If $X \nearrow Y$ is a conifold transition of Calabi--Yau threefolds, we have seen above (cf.\ \cite[Thm.\ 5.1]{HL12}) that 
\begin{align}\label{PTratioX}
    \psi_\ast \frac{\PTop (Y)}{\PTexc (Y)} = \PT (X)
\end{align}
by applying the variable change $\psi_\ast (\beta_Y, n) \coloneqq (\psi_\ast \beta_Y, n)$ to \eqref{PTratio}. It is natural to ask further whether the equality holds if $\oX$ has at worst terminal singularities.
\begin{qst}\label{qst}
Does the formula \eqref{PTratioX} hold for a small transition $X \nearrow Y$ of Calabi--Yau threefolds?
\end{qst}


\emph{Acknowledgments.} We thank Yukinobu Toda, Chin-Lung Wang, Baosen Wu and Zijun Zhou for helpful discussions. We also thank the referees for helpful comments and suggestions on this article. YL is supported by grants from the Fundamental Research Funds for the Central Universities and Applied Basic Research Programs of Science and Technology Commission Foundation of Shanghai Municipality (22JC1402700). SSW is supported by the National Science and Technology Council (NSTC) under grant number 111-2115-M-A49-019-MY3 and thanks the Institute of Mathematics at Academia Sinica for providing support and a stimulating environment.

\bibliography{GWPT}

\newcommand{\etalchar}[1]{$^{#1}$}

\end{document}